\documentclass[12pt]{article}
\usepackage{amsmath}
\usepackage{amsfonts}
\usepackage{amssymb}
\usepackage{theorem}

\title{Metric ultraproducts of finite groups with respect to some length functions}
\author{A.Ivanov 
\thanks{The research is supported by Polish National Science Centre grant DEC2011/01/B/ST1/01406} 
}

\newtheorem{theorem}{Theorem}[section]
\newtheorem{proposition}[theorem]{Proposition}

\newtheorem{lemma}[theorem]{Lemma}
\newtheorem{definition}[theorem]{Definition}

\newtheorem{remark}[theorem]{Remark}

\begin{document} 
\topmargin = 12pt
\textheight = 630pt 
\footskip = 39pt 

\maketitle

\begin{quote}
{\bf Abstract} 
We consider metric ultraproducts of finite groups 
with respect to some classes of length functions. 
All sofic groups embed into these ultraproducts. 
We study embeddings of normed groups. 
We also show that in some natural  situations 
such an ultraproduct is a simple group. 
\end{quote}

\section{Introduction}

Let $G$ be a group. 
A function $l:G\rightarrow [0,\infty)$ is called 
a {\em pseudo length function} if \\ 
(i) $l(1) = 0 $; \\ 
(ii) $l(g)=l(g^{-1})$; \\
(iii) $l(gh)\le l(g) + l(h)$. \\ 
A {\em length function} is a pseudo length function satisfying 
\begin{quote} 
$(i')$ $l(g) = 0 $ if and only if $g=1$, where $g\in G$.  
\end{quote} 
A pseudo length function is {\em invariant} if 
$l(h^{-1} gh) = l(g)$ for all $g,h\in G$. 
In this case it defines an invariant pseudometric by $l(gh^{-1})$. 
It becomes metric if $l$ is a length function. 
In this case we say that $G$ is a {\em normed group}. 
We consider normed groups as metric groups too. 

Metric ultraproducts of normed groups of 
bounded diameter, say $r$, are defined as follows.  
Let $(G_i ,l_i )$, $i\in I$, be a family of groups 
equipped with invariant length functions 
and let $\Phi$ be an ultrafilter on $I$. 
Then 
$$ 
N=\{ (x_i )_{i\in I}\in \prod_I G_i : 
lim_{i\rightarrow  \Phi} l(x_i ) = 0 \} 
$$ 
is a normal subgroup of $\prod_I G_i$.   
The {\em metric ultraproduct} $\prod_I (G_i ,l_i )/\Phi$ 
is defined to be $(\prod_I G_i )/N$ where 
the length function is defined by 
$$ 
l(xN) = lim_{i\rightarrow \Phi} l_i (x_i ) . 
$$ 
The latter is well-defined by compactness of $[0,r]$. 
This definition corresponds to Section 2.4 from \cite{Pestov}.

Metric ultraproducts of finite normed groups are  
deserved a particular attention in group theory.  
This is mainly motivated by investigations of 
{\em sofic groups}.  
We remind the reader that a group $G$ is called 
{\em sofic} if $G$ embeds into 
a metric ultraproduct of finite symmetric groups 
with the Hamming distance \cite{Pestov}. 
A group $G$ is called {\em hyperlinear} if $G$ embeds into 
a metric ultraproduct of finite-dimensional unitary  groups  
with the normalized Hilbert-Schmidt metric \cite{Pestov}. 
It is an open question whether these classes are the same 
and whether any countable group is sofic/hyperlinear. 

If in the definition of metric ultraproducts we 
do not assume diameter  boundedness, 
we arrive to a more general construction, 
see \cite{Pestov}, Sections 2.3 and  2.4, 
and to the following notion.    
A group $G$ is called {\em weakly sofic} if $G$ embeds 
into a metric ultraproduct of finite  groups with 
invariant length functions \cite{GR}. 
It is not known if this class coincides with the former ones. 

Thus the following general question becomes very interesting:  
{\em how strongly do properties of metric ultraproducts 
depend on particular choice of length functions?}   
In particular, {\em describe metric ultraproducts 
of familiar classes of finite groups with various 
invariant length functions. } 
 
It is worth noting that papers \cite{AP}, \cite{ST} 
and \cite{Thom} already suggest 
considering classes of finite groups with some 
special length functions. 
This in particular produces new versions of soficity. 

In our paper we concentrate on some natural 
modifications of the length function from the paper of 
A.Stolz and A.Thom \cite{ST}. 
On the one hand they are sufficiently  
general to cover all sofic groups. 
On the other hand some statements of \cite{ST} 
still hold with respect to them.  
Moreover we discover that in some respects these 
length functions resemble the Hamming length 
in symmetric groups. 
In Section 3 we prove that in the case 
of classical groups of unbounded dimension 
the corresponding metric 
ultraproducts are always simple groups. 
This generalises the corresponding theorem 
from \cite{ES} concerning symmetric groups 
with Hamming distance.  

When a metric ultraproduct $\prod (G_i ,l_i )/\Phi$ 
is a simple group there is no $\delta$ such that 
all elements of length $<\delta$ form a proper 
subgroup. 
Thus it becomes an interesting question if 
one can correct length functions $l_i$ so that 
the ultraproduct gets such a subgroup. 
Moreover given $G$, a finitely generated 
subgroup of  $\prod (G_i ,l_i )/\Phi$, 
{\em what norms can $G$ inherit under such 
embeddings?} 
We will discuss this in Sections 2 and 4.  
In Section 4  we  in particular show that 
under some mild assumptions on a subgroup 
$G_1 <G$ we can shorten the norms 
of elements of $G_1$ in such embeddings. 
Moreover we can do it using our modifications 
of length functions from \cite{ST}. 

We continue this introduction in Section 2.2 
after some definitions and easy observations 
of Section 2.1. 

The author is gtateful to Dugald Macpherson 
for interesting discussions concerning Section 3.

\section{Pseudo length functions} 

\subsection{Length corrections} 

The following lemma is obvious. 

\begin{lemma} 
Let $H<G$. 
For any invariant (pseudo) length function $l$ 
on $G$ the restriction of $l$ to $H$ is an invariant 
(pseudo) length function. 
\end{lemma}

The following statetments are taken from 
the paper of A.Stolz and A.Thom \cite{ST} (Lemmas 2.1 and 2.2).  

\begin{lemma} 
(i) If $G$ is a finite group with an invariant 
(pseudo) length function $l$ and $H$ is a normal 
subgroup of $G$, then  
$$ 
l_{G/H} (gH )= inf \{ l(gh): h\in H\} 
$$ 
defines an invariant (pseudo) length function on $G/H$. 

(ii) If $G$ is a finite group and $H$ is a normal subgroup of $G$
so that $G/H$ has an invariant pseudo length function $l$, then  
$$ 
l^{G} (g )= l(gH) 
$$ 
is an invariant pseudo length function on $G$. 
\end{lemma} 

%

Using this lemma we now describe some ways  
of correction of length functions.

\begin{lemma} 
Let $G$ be a finite group, $H$ be normal in $G$ and $l$ 
be an invariant (pseudo) length function bounded by $r$. 
Then \\ 
(i) for any natural $k\ge r$ and 
$s$ the following function $l^{*sk}_H$ 
is an invariant (pseudo) length function 
on $G$ bounded by 1: 
$$ 
\mbox{ if } g\in H \mbox{ then } 
l^{*sk}_H (g) := \frac{\sqrt[s]{l(g)}}{k + 1} 
\mbox{ and if } g\not\in H \mbox{ then } 
l^{*sk}_H (g) := 1 . 
$$  
(ii) if $r=1$, then for any natural numbers $s,s'$ and $t$ the function   
$$ 
l^{ss't}_{H} (g) = \frac{t}{t+1} \sqrt[s]{l_{G/H} (gH)} + \frac{\sqrt[s']{l(g)}}{t+1} 
$$ 
is an invariant (pseudo) length function bounded by $1$. 
\end{lemma}

{\em Proof.} 
The proof is straightforward; one should use the inequality 
$\sqrt[s]{a+b} \le \sqrt[s]{a} + \sqrt[s]{b}$. 
The boundedness follows from the inequality 
$ l^{ss't}_H \le max (\sqrt[s]{ l },\sqrt[s']{ l })$. 
$\Box$ 

\bigskip 

Let $\mathcal{G}_1$ and $\mathcal{G}_2$ be two 
countable families of groups with pseudo length functions. 
We assume that $\mathcal{G}_1$ and $\mathcal{G}_2$ 
are defined on the same class of groups, i.e. for any group $G$ 
$$ 
\mbox{there is } l_1 \mbox{ with } (G,l_1 )\in \mathcal{G}_1 
\Leftrightarrow 
\mbox{there is } l_2 \mbox{ with } (G,l_2 )\in \mathcal{G}_2 . 
$$ 
Let us enumerate all groups of this form: $G_0 ,G_1 ,G_2 ,...$ . 
We say that $\mathcal{G}_1$ is {\em asymptotically bounded}  
by $\mathcal{G}_2$ (with respect to our enumeration) 
if there are constants $c$ and $n_0$ 
so that for every $n>n_0$ and every choice of 
elements $g\in G_n$ we have $l_1 (g) \le c l_2 (g)$ 
(see \cite{ST}). 
We call these classes {\em asymptotically equivalent}  if 
they are asymptotically bounded with respect to each other. 

We say that $\mathcal{G}_1$ is {\em asymptotically bounded}  
by $\mathcal{G}_2$ {\em up to a polynomial}  
if there is a constant $c$ and natural numbers $m>0$ and $n_0$ 
so that for every $n>n_0$ and every choice of elements 
$g\in G_n$ we have 
$$ 
(l_1 (g))^{m} \le c \cdot l_2 (g) \mbox{ if } l_1 (g) < 1  \mbox{ and } 
$$  
 $$ 
(l_1 (g))  \le c\cdot l_2 (g)^{m}  \mbox{ if } l_1 (g) \ge 1 .  
$$  
The following lemma is obvious. 

\begin{lemma} \label{compare} 
Let $l\le 1$ be a generic pseudo length function on  
an enumerated family of all pairs of finite groups 
$(G,H)$ where $H$ is normal in $G$ 
(and $l$ is defined on $G$). 

Then for any fixed $t$ the pseudo length functions 
$l$ and $l^{11t}_H$ are asymptotically 
equivalent:  
$$ 
l \le (t+1) l^{11t}_H \mbox{ and } 
l^{11t}_H \le l . 
$$ 
On the other hand for any natural numbers $s,s'$ 
$$ 
l \le (t+1) l^{ss't}_H \mbox{ and } 
(l^{ss't}_H )^m \le l  \mbox{ where } m=max (s,s' ) ,  
$$  
i.e. the pseudo length functions 
$l$ and $l^{ss't}_H$ are asymptotically 
equivalent up to a polynomial. 
\end{lemma} 

The main object of our paper is the following 
pseudo length function. 

\begin{definition} 
If $G$ is a finite group then the function 
$$ 
\mbox{ {\bf conjugacy length} } 
l_c (g) = \frac{log |g^G |}{log |G|}
$$ 
is an invariant pseudo length function;  
it defines an invariant pseudometric by 
$l_c (g_1 g^{-1}_2 )$   
(Proposition 2.3 of \cite{ST}).
\end{definition} 
 
By Proposition 2.4 of \cite{ST} the function $l_c$ 
 is a length function exactly when $G$ has trivial centre. 

Let $H$ be a normal subgroup of $G$. 
Applying the lemmas above we obtain that 
\begin{quote} 
the functions $(l_c )^{*sk}_H$ and  
$(l_c )^{ss't}_{H}$ 
are invariant pseudo length function. 
\end{quote} 
We will see in Section 4 that they are very helpful. 

The same construction can be applied 
to some other pseudo length functions. 
The following invariant length function is taken from \cite{GR}. 
It will be used in the proof of Theorem \ref{WS} below.  
Consider a countable class of groups $\mathcal{C}$ 
where every group $G$ is considered as a pair 
together with a finite distingueshed set $\Delta$ 
of conjugacy classes. 
Then the conjugacy graph $\Gamma (G,\Delta )$ 
is defined on the set $V$ of all conjugacy classes 
of $G$ as follows: a pair $(x,y)\in V\times V$ 
is an edge if for some $c\in \Delta$ we have 
$x\subset cy$. 
We assume that for any $(G,\Delta )\in \mathcal{C}$ 
the diameter of the connected component 
containing $e_G$ is finite. 
\begin{quote} 
Define $l_{\Delta}(g)$ to be the distance 
from $e_G$ to $g^G$ in $\Gamma (G,\Delta )$. 
\end{quote} 
When $e_G$ and $g^G$ belong to distinct connected 
components of the conjugacy graph we assume that 
$l_{\Delta}(g)$ is the next natural number after 
the diameter of the component of $e_G$. 
This function can be normalised. 
If $G$ is finite and $H$ is the normal 
subgroup of $G$ generated by $\Delta$, then   
the definition of $l_{\Delta}$ slightly 
resembles $(l_c )^{*sk}_H$.    

\begin{remark} 
{\em 
Assume that $G=H$ and $S$ is a finite symmetric 
set representing all classes of $\Delta$ so that 
$G$ is generated by $S$. 
The {\em cancelation length} of a word $w$ in  
the alphabet $S$ is defined to be the least number 
of letters to be deleted from $w$ in order to obtain 
a word trivial in $G$. 
The {\em cancelation norm} of an element 
$g\in G$ is defined to be the minimal cancelation 
length of a representing word. 
By Proposition 2.A of \cite{BGKM} it coincides 
with $l_{\Delta}$. } 
\end{remark} 

\subsection{Subgroups of metric ultraproducts} 

Consider a metric ultraproduct 
$\prod_I (G_i ,l_i )/\Phi $ 
of finite normed groups. 
It is easy to see that replacing each 
$l_i$ by $\frac{l_i}{1+l_i}$ 
the ultraproduct becomes the same group 
with respect to a length function $< 1$. 
Below we will usually assume that 
length functions are bounded by $1$. 

Assuming that for every $i$, 
$H_i$ is a normal subgroup of $G_i$ 
we may replace 
$l_i$ by $(l_i )_{H_i}^{*11}$. 
As a result we obtain a new 
ultraproduct where all elements 
of length $<\frac{1}{2}$ 
form a normal subgroup. 
In particular it can be non-isomorphic 
with $\prod_I (G_i ,l_i )/\Phi $  
as an abstract group. 
For example this is the case for 
ultraproducts of classical groups 
with respect to the conjugacy length. 
This obviously follows from 
the main theorem of Section 3 
that when dimension is not bounded 
these ultraproducts are simple groups. 
Thus taking $G_i$ as $PGL_{n_i}(q_i )$ 
and taking $H_i$ as the corresponding 
$PSL_{n_i} (q_i )$ we obtain 
a nice illustration of the observation above. 

Let $P$ be a finitely generated 
group which embeds into
$\prod_I (G_i ,l_i )/\Phi $.  
It can happen that the length correction 
above forbids an embedding of 
$P$ into the metric ultraproduct 
after the correction.  

{\em How can we correct 
the length functions so that 
$P$ still a subgroup of 
the corresponding ultraproduct, 
but its length function 
becomes completely different?}  
In particular let $P_1$ be 
a distinguished normal subgroup 
of $P$. 
Can we arrange that $P_1$ 
consists of elements of 
length $<\frac{1}{2}$? 
This question will be studied in Section 4.

\section{Simplicity of metric ultraproducts} 

Let $P$ be a metric ultraproduct of 
normed groups of a class ${\cal C}$. 
If ${\cal C}$ is the class of all finite simple 
groups with the conjugacy length then by 
Theorem 3.1 of \cite{ST} the group $P$ is simple. 
Consider the situation when ${\cal C}$ consists 
of groups $G$ with trivial centre containing 
a subgroup $H$ which is simple. 
The best example is the class ${\cal S}$ of 
all finite symmetric groups $S_n$ with the subgroups 
$A_n$ of even permutations. 
It is shown in \cite{ES} that any infinite  
metric ultraproduct of groups from ${\cal S}$ 
with respect to Hamming's metrics is still 
a simple group. 
In this section we obtain the same result 
in some other cases. 
We start with the following logic lemma. 

\begin{lemma} \label{lemsim} 
Let ${\cal C}$ consist of finite normed 
groups $G$ with trivial centre 
and bounded norms.   
We assume that for every $\varepsilon >0$ 
there is a natural number $m$ such 
that any $G\in \mathcal{C}$ contains 
a simple subgroup $H$ which 
coincides with $(h^H )^m$ 
for any $h\in H$ with $l(h)>\varepsilon$. 

Assume that for any sequence 
$g_i \in G_i \in \mathcal{C}$ 
there is a sequence $h_i \in G_i$ 
such that $lim_{i\rightarrow \infty} l(h_i ) =0$, 
and all $g_i h_i$ belong to $H_i$. 

Then any infinite metric ultraproduct of 
$G_i$  is a simple group. 
\end{lemma} 

{\em Proof.} 
Let $\Phi$ be an ultafilter on $\omega$ 
such that the ultraproduct 
$\prod G_i /\Phi$ is infinite. 
For any sequence $g_i \in G_i \in \mathcal{C}$ 
such that under $\Phi$ the limit  
$lim_{i\rightarrow\infty} l_c (g_i )$ exists 
and is greater than $0$, there is 
a sequence $g'_i \in G_i$ such that 
$lim_{i\rightarrow \infty} l(g'_i ) =0$, 
and all $g_i g'_i$ belong to $H_i$. 
Thus any non-trivial element of 
the ultraproduct of $G_i$ can be 
presented by a sequence $h_i$ which 
belongs to the ultraproduct of 
the corresponding $H_i <G_i$.  

Assume that $g \in \prod G_i /\Phi $ 
is represented by $(g_i )_{i\in \omega}$ 
and $(h_i )_{i\in \omega}$ as above and 
$\varepsilon <lim_{i\rightarrow \Phi} l(g_i )$. 
Choose $m$ as in the formulation of the lemma. 
Since $(h^{H_i}_i )^m =H_i$ with 
respect to $\Phi$, we see that  
$$ 
(g^{\prod G_i /\Phi } )^m = \prod G_i /\Phi .   
$$  
$\Box$ 

\bigskip 

We now show how this lemma works in some very natural cases. 
We will apply it together with a very strong result from \cite{LS}. 
Theorem 1.1 of \cite{LS} states the existence of a universal 
constant $c$ such that  whenever  $H$ is a finite non-abelian 
simple group, $g\in H\setminus \{ 1\}$ and 
$m\ge c\frac{log |H|}{log |g^H |}$ 
we have $(g^H )^m = H$.   

Let $V$ be a an $n$-dimensional vector space over 
a finite field $K$ with $|K|=q$. 
Let $f$ be a sesquilinear form on $V$, write 
$G$ for the group of $K$-linear maps preserving $f$, 
write $Z$ for the centre of $G$ and write $H$ for 
the derived group of $G$. 
We are concerned with the cases when 
$f$ is one of the following: the zero form, 
a nondegenerate symplectic form, a nondegenerate hermitian form     
with respect to a non-identity involutory automorphism $J$ 
of $K$, or a nondegenerate symmetric form (for odd $q$). 
There is an additional case when $q$ is even, 
we then consider $G$ as a group of linear maps preserving a 
non-singular quadratic form $Q$ with $f(x,y) = Q(x+y) - Q(x)-Q(y)$.    
These are {\em classical groups} over finite fields. 

In all these cases the groups $HZ/Z$ are simple 
except for small values of $dim V$ and $|K|$.   
\begin{quote} 
Describe the metric ultraproducts of groups $G/Z$ 
with respect to natural invariant length functions which are simple.   
\end{quote} 
It is worth noting that by Proposition 2.7 of \cite{wilson} 
the ultraproducts of groups $G/Z$  with respect to 
the $\{ 0,1\}$-metric are not simple unless 
$G=H$. 
When the groups are normed by the conjugacy length 
the situation is different.   
The first statement of the following theorem 
is in fact proved in \cite{ES}.

\begin{theorem} \label{simple} 
Let $\mathcal{S}$ be the family of all 
symmetric groups $S_n$  and let $\mathcal{G}$ 
be the class of all projective classical groups $G/Z$ 
over finite fields of dimensions $>4$. 
We consider groups of both classes as normed 
groups with respect to the conjugacy length.  
Then \\ 
(i) any infinite metric ultraproduct of members 
of $\mathcal{S}$ is a simple group; \\
(ii) any non-discrete metric ultraproduct of 
members of $\mathcal{G}$ is a simple group; 
any discrete metric ultraproduct of members 
of $\mathcal{G}$ is a Chevalley  group over 
a standard ultraproduct of finite fields. 

The same statements hold for these classes with 
respect to length functions $(l_c )^{ss't}_H$ and 
the corresponding subgroups $A_n$ and $HZ/Z$ 
respectively (where $s$, $s'$ and $t$ are fixed).  
\end{theorem} 

{\em Proof.} (i) By Theorem 2.16 of \cite{ST} 
the conjugacy length is asymtotically equivalent 
to the Hamming's distance in $\mathcal{S}$ 
as well as in  $\mathcal{A}$, 
the family of finite alternating groups $A_n$.  
Now note that the Hamming's norm of any 
transposition $(i,j)$ is $\frac{2}{n}$. 
Thus the conjugacy length of a transposition from 
$S_n$ converges to $0$ when $n\rightarrow \infty$. 
This in particular means that for any 
sequence $g_i \in G_i \in \mathcal{S}$ 
such that the metric ultraproduct of $G_i$ 
is infinite and under the corresponding ultrafilter 
$lim_{i\rightarrow\infty} l_c (g_i )$ exists and 
is greater than $0$, there is a sequence $g'_i \in G_i$ 
such that $lim_{i\rightarrow \infty} l(g'_i ) =0$, 
and all $g_i g'_i$ are even permutations. 
Thus any element of the ultraproduct of $G_i$ 
can be presented by a sequence $h_i$ 
which belongs to the ultraproduct of 
the corresponding $A_n <S_n$. 
  
By Theorem 1.1 of \cite{LS} there is a universal constant 
$c$ such that for any $H\in \mathcal{A}$, $g\in H$ and 
$m\ge c\frac{log |H|}{log |g^H |}$ we have $(g^H )^m = H$.  
Since $|S_n |= 2|A_n |$ and $|g^{S_n}| \le 2|g^{A_n}|$, 
the equality  $(g^H )^m = H$ 
holds for $H=A_n$ under the assumption that 
$$ 
m\ge 2c\frac{log |S_n |}{log |g^{S_n} |} 
\mbox{ (which is } \frac{2c}{l_c (g)} ). 
$$ 
Thus to apply Lemma \ref{lemsim}, given $\varepsilon >0$ 
take $m\ge \frac{2c}{\varepsilon}$. 
Then the sequence $(g_i )$ as above defines 
an element of the metric ultraproduct such 
that the $m$-th power of 
its conjugacy class covers the group. 

By Lemma \ref{compare} this argument works 
for the length functions $(l_c )^{ss't}_{A_n}$.  \\ 

(ii) Consider a normed ultraproduct 
$\prod_{j\in J} G_j /\mathcal{U}$. 
Since $\mathcal{U}$ is an ultrafilter and 
$\mathcal{G}$ divides into finitely 
many classical series, we may assume  
that all $G_j$ belong to one of them. 

Let us start with the case 
of ultraproducts from $\mathcal{PGL}$, 
the class of all projective linear groups 
$PGL_n (q)$ over finite fields ${\bf F}_q$ 
with $n>4$. 
We show that any non-discrete metric 
ultraproduct of members of $\mathcal{PGL}$ 
is a simple group. 

By Theorem 2.20 of \cite{ST} in 
the class of all $PSL_n (q)$ and in the class 
of all $PGL_n (q)$ 
the conjugacy length is asymtotically 
equivalent to the {\em Jordan length} $l_J$, 
where 
$$ 
l_J (g) = n^{-1} inf_{a\in F^*_q} rk(a-g). 
$$ 
Note that the Jordan length does not 
depend on $q$, in particular the minimal norm 
of a non-trivial element from $PGL_n (q)$ 
is always $\frac{1}{n}$.  
This means that if the metric ultraproduct 
of groups $G_i \in \mathcal{PGL}$ 
is not a discrete space then the dimension 
of these groups is not bounded 
with respect to $\mathcal{U}$. 

Now note that the Jordan length 
of a diagonal $n\times n$-matrix 
$diag (1,1,...,1,d)$ 
is $\frac{1}{n}$. 
Thus the conjugacy length of such 
a matrix from $PGL_n (q)$ converges to $0$ 
when $n\rightarrow \infty$. 
This in particular means that for any 
sequence $g_i \in G_i \in \mathcal{PGL}$ 
such that the metric ultraproduct of $G_i$ 
is non-discrete and  
$lim_{i\rightarrow\mathcal{U}} l_c (g_i )$ 
is greater than $0$, there is 
a sequence of diagonal $h_i \in G_i$ such that 
$lim_{i\rightarrow \mathcal{U}} l(h_i ) =0$, 
and all $g_i h_i$ are from the corresponding $PSL_n (q )$-s. 
Since 
$\frac{log(q)}{log |PGL_n (q )|} \rightarrow 0$ 
with respect to $\mathcal{U}$ and   
$lim_{i\rightarrow\mathcal{U}} l_c (g_i h_i )$ 
is greater than $0$, we have  that 
$$ 
\frac{log (q )}{log |(g_i h_i )^{PSL_n (q )}|} \rightarrow 0  
\mbox{ under } \mathcal{U} . 
$$ 
By Theorem 1.1 of \cite{LS} there is a universal 
constant $c$ such that for any $H\in \mathcal{PSL}$, 
$g\in H$ and $m\ge c\frac{log |H|}{log |g^H |}$ 
we have $(g^H )^m = H$.  
Using the observations above, 
as in part (i) we obtain that the sequence 
$(g_i h_i )$ defines an element of the metric 
ultraproduct such that for some $m$ the $m$-th 
power of its conjugacy class covers the group. 
In fact to apply the argument of the proof of 
Lemma \ref{lemsim} it suffices to take 
$m\ge  \frac{2c}{l_c (g_i h_i )}$ with respect 
to $\mathcal{U}$ (computing $l_c$ in $G_i$),  i.e. 
$$ 
m \ge 2c\frac{log(q) + log |H |}{log (q ) + log |(g_i h_i )^{H} |} 
\mbox{ with respect to } \mathcal{U} . 
$$

If the metric ultraproduct of groups $G_i$ 
with respect to the Jordan length 
is infinite and discrete, then there is $n$ 
such that all $G_i$ are of the form $PGL_n$ 
up to the ultrafilter. 
As in this case the metric ultraproduct coincides 
with the non-metric one, we may apply 
the results of \cite{Point}.  
By Proposition 1 of that paper and some 
folklore observations, the ultaproduct of 
$PGL_n (K_i )$ (resp. $PSL_n (K_i )$) is 
$PGL_n$ ($PSL_n$) over the corresponding 
ultraproduct of fields. 

The use of the Jordan length above 
can be replaced by some 
additional computations. 
Such arguments will be applied below. 

Let us consider the cases of non-trivial 
bilinear forms (by Proposition 3.1 from \cite{ST} 
we may omit the case when $G/Z = HZ/Z$, i.e. for 
example where $f$ is symplectic over an odd field). 
By Section II.3  of \cite{D} 
the determinant of any unitary transformation 
over a field $K$ is of the form 
$\gamma^J \gamma^{-1}$ 
where $\gamma \in K$. 
This in particular means that multiplying any 
matrix from $U_n (K)$ by an appropriate 
diagonal matrix of the form  
$diag (1,1,...,1, \gamma^{J} \gamma^{-1})$ 
we obtain an element from $SU_n (K)$. 

The conjugacy length of such 
a diagonal matrix can be computed as follows. 
Writing the matrix as $e + (d-1)e_{nn}$ 
with $d= \gamma^{J} \gamma^{-1}$, 
we see that any its conjugate is of the 
form $e + (d-1)(\vec{c})'  \vec{b}$ where 
$\vec{b}$ and $\vec{c}$ are $n$-vectors 
so that $\vec{b} (\vec{c})' =1$.  
Since the unitary space over ${\bf F}_{q^2}$ 
always has an orthonormal basis and   
any matrix from $U_n (q^2 )$ takes 
an orthonormal bases to an orthonormal one, 
we have $\vec{b} = \vec{c}^J$. 
In particular the size of the conjugacy class 
$diag (1,1,...,1, \gamma^{J} \gamma^{-1})$  
in $U_n (q^2 )$
is between $q$  and  $q^{2n}$. 
Since by Section II.6 of \cite{D} 
$$
|U_n (q^2 )|= (q^n - (-1)^{n})q^{n-1}(q^n -(-1)^{n-1}) q^{n-2} 
\cdot ...\cdot (q^2 -1) q( q+1) , 
$$ 
the conjugacy length of 
$diag(1,...,1, \gamma^{J} \gamma^{-1})$ is 
asymptotically equivalent to $\frac{1}{n}$.  
This means that  if the metric ultraproduct 
of groups $G_i$ of the form $U_{n}$ 
is a discrete space then the dimension 
of these groups is bounded 
with respect to $\mathcal{U}$. 

A non-central element $g\in U_n$ with 
the minimal conjugacy length 
has the maximal size of its centraliser. 
We decompose $g$ into a product of 
commuting semisimple and unipotent elements 
and consider their Jordan forms. 
As in Section 7 of \cite{LS} we obtain 
a decomposition of the matrix of $g$ into 
a sum of diagonal matrices, Jordan blocks 
and tensor products of Jordan blocks with 
irreducible matrices.   
When the decomposition is non-trivial, 
the index of $C_{U_n} (g)$ in $U_n$ is 
clearly greater than $q$. 
Since by Lemma 4.4 of \cite{LS} the size of 
the centraliser of the Jordan $k$-block is not 
greater than $q^{2k}$, 
we also have  $|U_n : C_{U_n} (g)|\ge q$
when $g$ is a Jordan block or a tensor product. 
Thus the conjugacy length of $g$ 
does not depend on $q$, i.e. the minimal 
norm of a non-central element from $U_n (q^2 )$ 
converges to $0$ only when 
$n\rightarrow \infty$. 

This already means that if the metric ultraproduct 
of groups $G_i$ of the form $PU_{n}$ 
is not a discrete space then the dimension 
of these groups is not bounded 
with respect to $\mathcal{U}$. 
Thus we may apply the same argument as above 
(using an appropriate part of Theorem 1.1 of \cite{LS}). 

In the orthogonal case the determinant is 
$1$ or $-1$. 
In general the derived subgroup $\Omega_n$ is 
a proper subgroup of $O^+_n$ (of determinant $1$).  
The structure of 2-subspaces of $V$ becomes crucial 
in this case. 
A 2-subspace $P$ is called a {\em hyperbolic plane} if 
it has a basis $v_1$ and $v_2$ consisting 
of singular vectors (i.e. $Q(v_i )=0$) with 
$f(v_1 ,v_2) =1$. 
In the case of characteristic $\not=2$ 
(when $Q$ is determined uniquely by $f$) we will 
use the following fact from Section II.8 of \cite{D}. 
\begin{quote} 
Let $P$ be a hyperbolic plane in $V$. 
Then any element $g$ of $O_n$ can be presented 
as a product $sw$ with $w\in \Omega_n$ so that 
$s$ fixes the orthogal complement of $P$. 
\end{quote} 
Such an $s$ is determined by an orthogonal 
transformation of $P$. 
Thus the size of the conjugacy class of $s$ in 
$O_n (q)$ is bounded by $q^n q^{n-1}$. 
On the other hand by Section II.9 of \cite{D} for odd $n$, 
$$ 
|O^+_n (q )|= (q^{n-1} -1)q^{n-2}(q^{n-3} -1)q^{n-4} 
\cdot ...\cdot (q^2 -1) q , 
$$ 
and for $n=2m$, 
$$
|O^+_n (q)|= (q^{2m-1} - \varepsilon q^{m-1})(q^{2m-2}-1)q^{2m-3}  
\cdot ...\cdot (q^2 -1) q \mbox{ , with } \varepsilon \in \{ -1, 1\} .  
$$ 
Thus the conjugacy length of $s$   
in $O_n (q )$ converges to $0$ when $n\rightarrow \infty$. 
In particular, the arguments above also work in this case 
(using appropriate places of \cite{LS} and \cite{Point}).

In the case of characteristic 2 an orthogonal space 
is determined by a quadratic form. 
In this case $f(x,y)$ is simplectic of even dimension. 
Moreover there are two kinds of quadratic geometry on $V$: 
an $O^{+1}$-geometry, where $V$ is an orthogonal sum of 
hyperbolic planes, and an $O^{-1}$-geometry, where $V$ 
is an orthogonal sum of $\frac{n}{2} -1$ hyperbolic planes 
and one {\em definite} (non-singular) plane.  
In both cases the statement on decomposition 
$g=sw$, with $w\in \Omega^{\varepsilon}$ also holds 
(\cite{D}, Section II.10). 
Since in the first case 
$$ 
|O^+_{2m} (q )|= (q^{m} -1)(q^{2(m-1)} -1)q^{2(m-1)}  
\cdot ...\cdot (q^2 -1) q^2 , 
$$ 
and in the second one 
$$ 
|O^+_{2m} (q )|= (q^{m} +1)(q^{2(m-1)} -1)q^{2(m-1)}  
\cdot ...\cdot (q^2 -1) q^2 , 
$$ 
we can apply our arguments above. 

By Remark \ref{compare} the arguments above also work  
for the length functions of the fom $(l_c )^{ss't}_{PSL_n}$ 
and their relatives. 
$\Box$

\bigskip

\begin{remark} 
{\em It is worth noting that if we consider 
the classes $\mathcal{S}$ and $\mathcal{PGL}$ 
with respect to the family of all length functions 
$(l_c )^{ss't}_H$ for all natural $s,s'$ and $t$ 
(with respect to $\mathcal{A}$ and $\mathcal{PSL}$), 
then Theorem \ref{simple}  does not remain true. 
For example if we consider $S_n$ with respect to 
$(l_c )^{n,1,3}_{A_n}$, then the norm of 
a transposition $(i,j)$  can be evaluated as 
$\frac{3\sqrt[n]{2d}}{4\sqrt[n]{n}} + \frac{d}{2n}$ 
for some constant $d$ (replace the conjugacy length by 
Hamming length and apply Lemma \ref{compare}). 
In particular for sufficiently large $n$ the norm of 
any element of $S_n \setminus A_n$ is greater than 
$\frac{1}{2}$. 
On the other hand by the definition of 
$(l_c )^{n,1,3}_{A_n}$ the norm of any element of $A_n$ 
is less than $\frac{1}{4}$. 
This in particular shows that the elements of 
the metric ultraproduct of $(S_n ,(l_c )^{n,1,3}_{A_n})$ 
which have norm $<\frac{1}{4}$, form a non-trivial 
normal subgroup. 
}
\end{remark}

\section{Embeddings of weakly sofic groups and LEF-groups}

Let $P$ be a finitely generated 
subgroup of a metric 
ultraproduct $\prod_I (G_i ,l_i )/\Phi $ 
of finite normed groups 
(i.e. $P$ is weakly sofic). 
Let $P_1$ be a distinguished 
normal subgroup of $P$ so that 
$P/P_1$ is residually finite. 
Using the approach of Section 2 
we now show  that $P$ 
can be embedded into a metric 
ultraproduct $\hat{P}$ 
of finite normed groups so that 
$P_1$ consists of all 
elements of $P$ which have 
norms $\le \frac{1}{2}$.  
Moreover we will show that 
when $P$ is a {\em LEF}-group  
\cite{GV}, the norms of 
that finite groups can be taken 
in the form $(l_c )^{*st}$, 
introduced in Section 2.  

For convenience of the reader we 
recall the following definition.  

\begin{definition} 
Let ${\cal C}$ be a class of normed groups. 
A group $G$ is said to have the ${\cal C}$-{\bf approximation property} 
\cite{Thom} if   for any $g\in G$ there exists $\delta_g>0$ 
such that for all finite $D\subset G$ and $\varepsilon >0$ there exists 
a group $(C,l)\in \mathcal{C}$ and a map $\phi :G\rightarrow C$ 
so that 
$$
\phi(e)=e \mbox { , }  l(\phi (g))\ge \delta_g \mbox{ for all } g\in D \mbox{ and } 
$$ 
$$ 
l(\phi (gh)\phi (h)^{-1}\phi (g)^{-1}) <\varepsilon \mbox{ for all } 
g,h \in D \mbox{ with } gh\in D . 
$$
\end{definition} 
Proposition 1.8 of \cite{Thom} states that 
a countable group has the ${\cal C}$-approximation property 
if and only if $G$ embeds ito a metric ultraproduct 
of members of $\mathcal{C}$ with respect to 
a non-principal ultrafilter on $\omega$.

\begin{theorem} \label{WS} 
Let $N\lhd N_1 \lhd F$, where $F$ 
is a finitely generated free group.  
Let $P=F/N$ and $P_1 =N_1 /N$. 
Assume that $F/N$ is weakly sofic 
and $F/N_1$ is residually finite. 
Then $P$ embeds into a metric ultraproduct 
of finite groups with norms $\le 1$ so that 
$P_1$ consists of all elements of $P$ which have 
norms $\le \frac{1}{2}$ in that ultraproduct.  
\end{theorem} 

{\em Proof.} 
The proof is based on results 
and methods of \cite{GR}. 
It also uses the idea of 
length corrections of Section 2. 
Let $\varepsilon <\frac{1}{4}$. 
We want to apply an appropriate 
version of the approximation property 
for $P$ and $\frac{10}{9} \varepsilon$ 
where $\delta_g = \frac{9}{40}$.  
We additionally arrange that 
all elements of $P_1$ have norms 
$\le \frac{1}{2}$ and elements of 
$P\setminus P_1$ are of norms 
$\ge \frac{21}{40}$.  

The resulting $(C,l)$ from the formulation 
will be denoted $(G,\tilde{l})$ below.  
Let $D$ be a finite subset of $F/N$ presented 
by words $\{ w_1 ,...,w_k \}\subset F$. 
We may assume that the empty word belongs to $D$.   
Let $r= 3 max (|w_i | :i\le k)$. 

Since $P$ is weakly sofic, by Theorem 4.3 
of \cite{GR} for any $u_1 ,...,u_m$ 
from $N$ the profinite closure of 
$u^F_1 \cdot ...\cdot u^F_m$ 
is contained in $N$. 
Thus there is a homomorphism $\phi$ from 
$F$ to a finite group $G$ so that for 
any $w\in D^3 \setminus N$ the element 
$\phi (w)$ does not belong to 
any product 
$\phi (u_1 )^G \cdot ... \cdot \phi (u_m )^G$ 
with $u_i \in N\cap D$ and 
$m\le \frac{9}{10\varepsilon}$. 
Let $H$ be a normal subgroup of $G$ 
defined by $H=\phi (N_1 )$. 
Since $F/N_1$ is residually finite  
we may suppose that 
$\phi (D\setminus N_1 )\cap H =\emptyset$.

Repeating the sufficiency argument 
of Lemma 6.4 of \cite{GR} we build 
an invariant  length function $l$ on $G$
so that all elements of $\phi (D^3 \cap N)$ 
have norm $\le \varepsilon$ but 
all elements of $\phi (D \setminus N)$ 
have norm $\ge \frac{9}{10}$. 
For convenience of the reader we 
recall that this is $\varepsilon l_{\Delta}$ 
(see Section 2) where $\Delta$ consists 
of conjugacy classes $\phi (w)^G$ 
with $w\in N$ and $|w|\le r$.

Let us define an invariant length 
function $\tilde{l}$ on $G$. 
We start with a preliminary 
function $\gamma (x)$ on $\mathbb{R}^+$. 
Let $m_0 = max (l(g):g\in G)$. 
Let $s>1$ satisfy $\sqrt[s]{m_0}<2$. 
Let $\gamma (x)$ be a continuous function 
which equals $\sqrt[s]{x}$ for real 
numbers $>\frac{9}{10}$ and is of the 
form $p \cdot x$ for $x\le \frac{9}{10}$. 
It is clear that $1\le p\le \frac{10}{9}$. 
It is worth noting that 
$\gamma (x +y)\le \gamma (x) + \gamma (y)$. 

For elements $g\in H$ we define 
$$
\tilde{l} (g) =\frac{1}{4} \gamma (l(g)). 
$$ 
When $g\not\in H$ let 
$$
\tilde{l} (g) =\frac{1}{4} \gamma (l(g)) + 
\frac{1}{3} \gamma (inf (l(gh):h\in H)). 
$$ 
By the choice of $l$, $H$ and $\gamma$ 
the function $\tilde{l}$ 
is an invariant length function. 

It is now clear that for any 
$w\in D \cap N_1$ the value 
$\tilde{l} (\phi (w))$ 
is less than $\frac{1}{2}$. 
Moreover when $w\in D^3 \cap N$ then 
$\tilde{l} (\phi (w))< \frac{10}{9} \varepsilon$, 
and when $w\in D\setminus N$ then 
$\tilde{l} (\phi (w))\ge \frac{9}{40}$. 

Note that when $g\in \phi (D)\setminus (H\cup N)$,    
the value $\tilde{l} (g)$ is greater than 
$\frac{9}{40} + \frac{3}{10}=\frac{21}{40}$, 
i.e. $> \frac{1}{2}$. 

On the other hand if 
$F/N \models (wN = vN \cdot uN )$, 
where $w,v,u \in D$, then 
$w^{-1}vu\in N\cap D^3$, 
i.e. the corresponding distance between 
$\phi (w)$ and $\phi (v) \phi (u)$ 
is not greater than $\frac{10}{9} \varepsilon$. 
$\Box$

\bigskip

We now describe some 
natural situations where functions  
$(l_c )^{*st}_H$ and $(l_c )^{ss't}_H$ appear.  
By Theorem 2.16 of \cite{ST} the conjugacy 
length $l_c$ and the Hamming length  
are asymptotically equivalent in the class 
$\mathcal{S}$ of all symmetric groups $S_n$. 
This in particular implies that any sofic 
group embeds into a metric ultraproduct of 
normed groups $(S_n ,l_c )$. 
By Lemma \ref{compare} $l_c$ can be replaced 
by any $(l_c )^{*st}$ or $(l_c )^{ss't}_H$ 
for fixed $s,s',t$ 
and $H\in \{ A_n ,S_n \}$. 

On the other hand functions of this 
kind can be applied for embeddings 
of sofic groups into some special 
metric ultraproducts. 
We now demonstrate it in the case 
of {\em LEF groups}. 

\begin{definition} 
(\cite{GV})
A group $G$ is called LEF if for every 
finite subset $D\subseteq G$ there is 
a finite group $C$ and an injective map 
$\phi : D\rightarrow C$ such that 
any triple $h,g,hg\in D$ satisfies 
$\phi (hg )=\phi (h) \phi (g)$. 
\end{definition} 

We think that it is folklore that  LEF 
is equivalente to statement (i) of 
the following proposition.  
We give the proof below just for 
a curious observation that LEF is 
equivalent to (ii) too. 
This is based on Theorem 4.3 of  
\cite{Ivanov}. 

\begin{proposition} \label{LEF1} 
Let $P$ be a finitely generated group. 
The group $P$ is LEF if and only if 
one of the following equivalent conditions holds: \\ 
(i) any presentation $P=F/N$ so that 
$N$ is a normal subgroup of a finitely 
generated free group $F$, 
satisfies the following property:  
\begin{quote} 
for any two finite subsets $D_1 \subseteq N$ 
and $D_2 \subseteq F\setminus N$ there 
exists a normal subgroup $H<F$ of finite index 
so that $D_1 \subseteq H$ and $D_2 \cap H=\emptyset$; 
\end{quote}  
(ii) the (coloured) Cayley graph of $P$ has 
the finite model property: 
any sentence of the first order theory of 
this graph holds in a finite graph. 
\end{proposition} 

{\em Proof.} By Theorem 4.3 of \cite{Ivanov} 
the conditions (i) and (ii) are equivalent. 
The condition (i) obviously implies LEF. 
To see that LEF implies (ii) we use the following 
consequence of Proposition 1.1  of  \cite{Ivanov}:  
\begin{quote} 
the Cayley graph $\Gamma_P$ of a finitely 
generated group $P$ has the finite model property 
if and only if for every natural $n$ there is a 
finite graph $\Gamma_n$ so that 
any point $v\in \Gamma_P$ and 
any $v'\in \Gamma_n$  have isomorphic 
$n$-balls. 
\end{quote}  
Now note that when $P$ is LEF 
then choosing the $n$-ball $B_n (1)$ of  
the neutral element $1$ as a finite subset  
$D\subseteq P$ from the definition of LEF
we find a finite group containing 
the partial group $B_n (1)$. 
Taking the subgroup generated by $B_n (1)$ 
if necessary, we obtain a finite group 
such that its Cayley graph has $n$-balls 
isomorphic to $B_n (1)$. 
We see that $\Gamma_P$ has 
the finite model property. 
$\Box$ 

\bigskip

Let $P$ be a finitely generated 
group which is LEF.  
Let $P_1$ be a distinguished 
normal subgroup of $P$, 
say $P= F/N$ and $P_1 = N_1 /N$, 
where $F$ is a finitely generated 
free group and $N\lhd N_1 \lhd F$. 

\begin{definition} 
A subgroup $P_1 <P$ is called 
LEF-separated if for any 
finite subsets $D_1\subseteq N$ 
and $D_2 \subset F\setminus N$ 
there is a homomorphism $\phi$ 
from  $F$ to a finite group $C$ 
so that $\phi (D_1 )=\{ e\}$, 
$\phi$ is injective on $D_2$ and 
$$ 
\phi (D_2 \setminus N_1 ) \cap \phi (N_1 ) =\emptyset .
$$  
\end{definition} 
 
It is clear that if $P_1$ is of finite index, 
then $P_1$ is LEF-separated. 
 
The following observation shows 
that when $P$ is centreless and 
$P_1$ is LEF-separated we can 
embed $P$ into a metric 
ultraproduct $\hat{P}$ of 
normed groups introduced 
in Section 2 so that $P_1$ 
consists of all elements of $P$ 
which have norms $\le \frac{1}{2}$.

\begin{proposition} \label{LEF2} 
Let $N\lhd N_1 \lhd F$, where $F$ 
is a finitely generated free group.  
Let $P=F/N$ be LEF, centerless and 
$P_1 =N_1 /N$ be LEF-separated in $P$. 

Then $P$ embeds into a metric 
ultraproduct of finite groups 
with trivial centre and with metrics 
defined by pseudo length functions 
of the form $(l_c )^{*st}_{H}$ 
so that $P_1$ 
consists of all elements of $P$ 
which have norms $\le \frac{1}{2}$ . 
\end{proposition} 

{\em Proof.}    
We want to verify the appropriate version 
of the approximation property for $P$ where   
$\varepsilon =0$, $\delta_g = \frac{1}{4}$ 
and the length of elements from $P_1$ 
is $\le \frac{1}{2}$.  
Let $D$ be a finite subset of $F/N$ presented 
by words $\{ w_1 ,...,w_k \}\subset F$. 
We may assume that the empty 
word belongs to $D$.   
For any $w\in D^3$ which does not 
present the neutral element find $w' \in F$ 
with non-trivial $[w,w']$ in $F/N$. 
Let $D'$ consist of all words $w\in D^3$  
and the corresponding $[w,w']$ 
of that form. 

By LEF-separation there is a homomorphism 
$\phi$ from $F$ to a finite group $G$ 
so that $D'\cap N \subset Ker \phi$,  
$(D' \setminus N) \cap Ker \phi =\emptyset$  
and 
$\phi (D'\setminus N_1 ) \cap \phi (N_1 ) =\emptyset$. 
Let $H =\phi (N_1 )$. 

Note that for any 
$h\in H$ the value 
$(l_c )^{*s1}_{H} (h)$ 
is less than $\frac{1}{2}$. 
Since for $w\in D^3 \setminus N$, 
the element $\phi (w)$ does not 
belong to the center of $G$,  
for any non-trivial $g\in \phi (D^3)$ 
the value $|g^G|$ is greater than $1$. 
Thus we may choose $s$ so that 
$(l_c )^{*s1}_{H} (g)> \frac{1}{4}$ 
for any non-trivial $g\in \phi (D^3 )$. 
Thus 
$(l_c )^{*s1}_{H} (g)> \frac{1}{4}$ 
for all $g\in G$ and  
$(l_c )^{*s1}_{H} (g)> \frac{1}{2}$ 
for $g\not\in H$. 

On the other hand if 
$F/N \models (wN = vN \cdot uN )$, 
where $w,v,u \in D$, then 
$w^{-1}vu\in N\cap D'$, 
i.e. $\phi (w) =\phi (v) \phi (u)$. 

$\Box$ 

\bigskip 

\begin{remark} 
{\em 
It is worth noting that 
if in the proof above 
when we apply LEF-separation 
we can additionally 
obtain that 
for any  $g\in \phi (D)$ 
the coset $gH$ does not intersect 
$Z(G)$, we can replace 
functions $(l_c )^{*s1}$ by 
$(l_c )^{ss't}$ for appropriate  
$s$ and $s'$. 
Indeed, as above we may choose $s'$ so that 
$(l_c )^{ss'1}_{H} (g)> \frac{1}{4}$ 
for any non-trivial $g\in \phi (D^3 )$. 
Moreover we can now choose $s$ so that 
$(l_c )^{ss'1}_{H} (g)> \frac{1}{2}$ 
for $g\in \phi (D^3 ) \setminus H$.  }
\end{remark}

\bigskip

University of Wroc{\l}aw, Poland 

E-mail: ivanov@math.uni.wroc.pl


\begin{thebibliography}{99}  
\bibitem{AP} G.Arzhantseva and L. Paunescu, 
Linear sofic groups and algebras, arXiv: 1212.6780.  
\bibitem{BGKM} M.Brandenbursky, \'{S}.R.Gal, J.Kendra, M.Marcinkowski, 
Cancelation norm and the geometry of biinvariant word metrics. 
ArXiv: 1310.2921.  
\bibitem{D} J.Dieudonn\'{e}, La g\'{e}ometrie des groupes classiques, 
(Springer, Berlin-G\'{o}ttingen-Heidelberg, 1955)
\bibitem{ES} G.Elek and E.Szabo, Hyperlinearity, 
eventually free actions and $L^2$-invariants, 
Math. Ann. 332(2005), 421 - 441. 
\bibitem{GV} E.Gordon and A.Vershik, Groups that are loally 
embeddable in the class of finite groups, 
Algebra i Analiz 9(1997), no. 1 71 - 97. (Russian)  
\bibitem{GR} L.Glebsky and L.M.Rivera, Sofic groups and 
profinite topology on free groups, J.Algebra 320(2008), 3512 - 3518.   
\bibitem{Ivanov} A.Ivanov, Cayley graphs having nice enumerations, 
Israel J. Math., 137(2003), 61 - 108. 
\bibitem{LS} M.Liebeck ana A.Shalev, Diameters of finite 
simple groups: sharp bounds and applications, Ann. Math.  
154(2) (2001), 384 - 406. 
\bibitem{Pestov} V.Pestov, Hyperlinear and sofic groups: a brief guide, 
Bull. Symb. Logic, 14(2008), 449 - 480. 
\bibitem{Point} F.Point, Ultraproducts and Chevalley groups, 
Arch. Math. Log., 38(1999), 335 - 372. 
\bibitem{ST} A.Stolz and A.Thom, On the lattice of normal 
subgroups in ultraproducts of compact simple groups, arXiv: 1207.0977. 
\bibitem{Thom} A.Thom, About the metric approximation of Higman's group, 
J.Group Theory, 15(2012), 301 - 310. 
\bibitem{wilson} J.Wilson, On simple pseudofinite groups, J. London Math. Soc. (2), 
51 (1995), 471 - 490. 
\end{thebibliography}
\end{document}